\documentclass[letterpaper, 10 pt, conference]{ieeeconf}  

\IEEEoverridecommandlockouts                              

\usepackage[utf8]{inputenc}

\usepackage{graphicx}
\usepackage{amsmath}
\usepackage[version=4]{mhchem}
\usepackage{siunitx}
\usepackage{longtable,tabularx}
\usepackage{cite}
\usepackage{amsmath,amsfonts,derivative}
\usepackage{subcaption}
\usepackage{algorithm}
\usepackage{algorithmic}

\usepackage{graphicx}
\usepackage{textcomp}
\usepackage{xcolor}
\usepackage{comment}
\def\BibTeX{{\rm B\kern-.05em{\sc i\kern-.025em b}\kern-.08em
    T\kern-.1667em\lower.7ex\hbox{E}\kern-.125emX}}

\usepackage{pgfplots}
\pgfplotsset{compat=1.18}
\usepackage{graphicx} 
\usepackage{tikz}
\usetikzlibrary{calc, patterns, decorations.pathmorphing,decorations.markings}
\usepackage{algorithm}

\title{\LARGE \bf
Switched Optimal Control with Dwell Time Constraints
}

\author{Masoud S. Sakha$^{1}$ and Rushikesh Kamalapurkar$^{2}$
\thanks{*This research is supporting, in part, by the Office of Naval Research grant N00014-21-1-2481.Any opinions, findings, conclusions, or recommendations
detailed in this article are those of the author(s), and do not necessarily
reflect the views of the sponsoring agencies.}
\thanks{$^{1}$Masoud S. Sakha and  Rushikesh Kamalapurkar are  with the Department of Mechanical and Aerospace Engineering, University of Florida, Gainesville, Fl, USA.
        {\tt\small masoud.sakha@ufl.edu, rkamalapurkar@ufl.edu}}.%
}

\begin{document}

\maketitle
\thispagestyle{empty}
\pagestyle{empty}

\begin{abstract}

This paper presents an embedding-based approach for solving switched optimal control problems (SOCPs) with dwell time constraints.  
At first, an embedded optimal control problem (EOCP) is defined by replacing the discrete switching signal with a continuous embedded variable that can take intermediate values between the discrete modes. While embedding enables solutions of SOCPs via conventional techniques, optimal solutions of EOCPs often involve nonexistent modes and thus may not be feasible for the SOCP. 
In the modified EOCP (MEOCP), a concave function is added to the cost function to enforce a bang-bang solution in the embedded variable, which results in feasible solutions for the SOCP. However, the MEOCP cannot guarantee the satisfaction of dwell-time constraints.

In this paper, a MEOCP is combined with  a filter layer to remove switching times that violate the dwell time constraint. Insertion gradients are used to minimize the effect of the filter on the optimal cost.
\end{abstract}

\section{Introduction}
In recent decades, switched systems have received increasing attention from scholars across various disciplines due to their capacity to model various real-world phenomena \cite{SCC.Shaker.Wisniewski2009, SCC.BechBorchersen.Larsen.ea2015, SCC.Seyed Sakha.Shaker.ea2018, SCC.Hespanha.Morse1999, SCC.Liberzon.Morse1999, SCC.Shorten.Wirth.ea2007, SCC.Lin.Antsaklis2009, SCC.Zhao.Zhang.ea2012, SCC.Leth.Wisniewski2012, SCC.Sloth.Wisniewski2014, SCC.Greene.Sakha.eatoappear}.
 Switched systems are a subset of hybrid dynamical systems, comprising a finite collection of subsystems and a logical framework governing transitions between them. Due to the presence of both discrete and continuous variables, switched systems can exhibit complex dynamic behaviors. 

Optimal control problem of switched systems involves the determination of the best control input, including the switching logic, to optimize performance criteria while accounting for discrete transitions among different subsystems (or modes) of the system. Due to the hybrid dynamical behavior of switched systems, conventional optimal control methods, based solely on the gradient of the cost function, cannot be used to solve such problems ~\cite{SCC.Xu.Antsaklis2003, SCC.Zhu.Antsaklis2014, SCC.Zhao.Gan.ea2020}.
 
Typically, switched systems can be categorized into two primary groups: those with externally forced switching (EFS), the focus of this study, and those with internally forced switching (IFS) \cite{SCC.Zhu.Antsaklis2014}. Methods that address switched optimal control problems (SOCPs) with EFS fall into two main categories: two-stage optimization and embedding.

Two-stage optimization methods, which were first introduced in \cite{SCC.Xu.Antsaklis2000}, use two layers of optimization where the inner layer optimizes the switching time instants and the outer layer optimizes the switching sequences. Many of the works in this category just focus on the first layer and optimize the switching times with respect to fixed switching sequences \cite{SCC.Xu.Antsaklis2000, SCC.Xu.Antsaklis2003, SCC.Xu.Antsaklis2004}. The main challenge for optimization of switching sequences is that the number of switching sequences grows combinatorially with respect to the number of switching events. A  mode insertion technique has been developed in \cite{SCC.Axelsson.Wardi.ea2005}  to avoid such combinatorial explosion, where in each interaction, the switching sequence is updated by inserting a new mode at a suitable time instance and for a suitable time period. In \cite{SCC.Axelsson.Wardi.ea2007}, a gradient descent approach has been used for mode insertion, and two main theoretical questions have been discussed. First, the meaning of convergence in Euclidean space of increasing dimensions, and second, the existence of algorithms that converge in that sense. The idea is extended to insert multiple new modes in each iteration in \cite{SCC.Wardi.Egerstedt.ea2014}. A new framework for convergence analysis is also introduced in \cite{SCC.Wardi.Egerstedt.ea2014} based on Polak’s notions of optimality developed in \cite{SCC.Polak.Wardi1984}   for infinite-dimensional optimization problems.


In contrast to two-stage optimization, embedding approaches treat the discrete mode sequence variable as a collection of continuous decision variables.
The idea of embedding switched systems into a larger family of continuous systems was initially presented in \cite{SCC.Bengea.DeCarlo2005}, which showed that if a bang-bang solution of the embedded optimal control problem (EOCP) exists, then it can be directly identified as a solution of the original SOCP. If not, a chattering lemma can provide a sub-optimal solution for the SOCP.
The advantage of the embedding approach over two-stage optimization is that there is no need to make any assumptions about the number of switching events.  It also enables the use of conventional optimal control methods like Pontryagin's minimum principle or dynamic programming.
Several numerical algorithms have also been explored in \cite{SCC.Wei.Uthaichana.ea2007} to solve embedded non-autonomous EOCPs, and it has been shown that sequential quadratic programming (SQP) based methods for solving EOCPs have a better performance than mixed integer programming (MIP) base methods to solve EOCPs.
However, a solution of the EOCP is not necessarily feasible for the SOCP since switching signals only accept their designated integer values, while the EOCP allows any real values in predefined intervals. 
To address the feasibility issue, \cite{SCC.Abudia.Harlan.ea2020} adapted the approach from \cite{SCC.Bengea.DeCarlo2005} by adding a concave auxiliary function to the cost function. 
This modification ensures that the modified EOCP (MEOCP) possesses a bang-bang solution, making it feasible for the SOCP.
%


Despite extensive research in SOCPs, there are still several open questions. One of the primary concerns in practical applications is the dwell-time constraint on the switching rate, which arises from various factors, including actuator rate limits.
For instance, high-speed switching in electrical actuators like relays can cause electromagnetic interference or overheating, while actuators with mechanical components, such as servo motors, can experience wear and tear, thereby reducing their lifespan. 
For example, short cycling in heating, ventilation, and air Conditioning (HVAC) systems can harm the system and decrease its efficiency, leading to temperature fluctuation and poor dehumidification. The dwell-time constraint helps relieve these issues by limiting the switching frequency.

 In \cite{SCC.Abudia.Harlan.ea2020}, it has been shown that dwell-time constraints can be heuristically met by adjusting the constants used to solve the MEOCP within the auxiliary cost function.
 
 Nonetheless, this method is heuristic, relies heavily on the numerical optimization algorithm, and cannot guarantee satisfaction of the dwell-time constraint. In this paper, instead of tuning the constants,  an extra layer has been added to the MEOCP solution to guarantee the satisfaction of the dwell time constraint. This layer removes any switching signals that violate the dwell time constraint in an approximately optimal manner. Simulation results are provided to illustrate the efficacy of the proposed approach. 
 
\section{Switched Optimal Control Formulation}
Consider the switched nonlinear system with two subsystems defined as:
\begin{equation}
		\dot{x} =f_{\bar{v}}(t,x,u)
   \label{Eq_switching function}
\end{equation}
Where $t \in [t_0,t_f]$ represents time, $x \in \mathbb{R}^{n}$ is a state vector, 
   	$\bar{v}: \mathbb{R}_{\geq0}\rightarrow  \bar{\mathcal{V}}:= \{ 0,1\}$ is a switching signal, $u \in \mathbb{R}^{m}$ is the control input, and $f_1, f_2 \in \mathcal{C}^1(\mathbb{R}^{n},\mathbb{R}^{n})$.

\vspace{2mm}
\textbf{Definition 1:} A switching signal  $\bar{v}: [t_0,t_f] \rightarrow \bar{\mathcal{V}}$ is called feasible if it is piecewise constant.
Let $(v^1,v^2, \cdots ,v^{N+1} )$, for some $N\geq 0$, denote the successive values of a feasible switching signal $\bar{v}$ and  let $\tau_i$ denote the time of the switch from $v^i$ to  $v^{i+1}$, such that $t_0 \leq \tau_1 \leq \cdots  \leq \tau_N \leq t_f$.

\vspace{2mm}
\textbf{Definition 2:} A feasible schedule for the switched system (\ref{Eq_switching function}) is denoted by $\sigma=(q, \tau)$ where $\tau = (\tau_1, \cdots,\tau_N)^{\top}$ is a vector of  switching time instances and  $q \in Q$, where $Q$ is a set of all sequences of elements from $\bar{\mathcal{V}}$ of size of $N+1$.

\vspace{2mm}
A SOCP with a dwell time constraint is formulated as 
 \begin{equation*}
	\begin{aligned}
		\underset{u(\cdot),\bar{v}(\cdot)}{\min} ~~& J(t_0,x_0,t_f,u(\cdot),\bar{v}(\cdot))  \\
		\text{subject to}: & \\
		&  \text{(i) $(t_0,x_0,t_f,x(t_f)) \in \mathcal{B}\subseteq \mathbb{R}^{2n+2} $ },\\
		&  \text{(ii) $\bar{v}(t)\in  \bar {\mathcal{V}} , u(t)\in  \mathcal{U} \in \mathbb{R}^m, ~ \forall t \in[t_0,t_f]$ }, \text{and}\\
		&  \text{(iii) $\forall t_1, t_2 \in[t_0,t_f]$ with $\bar{v}(t_1^-)\neq \bar{v}(t_1^+)$, and}~\\
		&~~~~~\bar{v}(t_2^-)\neq \bar{v}(t_2^+) , |t_1-t_2|\geq T>0 ,
	\end{aligned}
\end{equation*} 
where (iii) encodes the dwell-time constraint, $\mathcal{B}$ and $\mathcal{U}$ are compact sets, and the cost function is defined as
 \begin{equation}
	\begin{aligned}
 J(t_0,x_0,t_f,u(\cdot),\bar{v}(\cdot)):= & \int_{ t_0   }^{ t_f  }  l_ {\bar{v}(t)}(t,x(t),u(t)) \mathrm{d}t \\
 & + K(t_0,x_0,t_f,x(t_f)),
	\end{aligned}
 \label{Eq_Cost SW}
\end{equation} 
where $l_{\Bar{v}(t)}$ is defined as
\begin{equation}
	\begin{aligned}
l_{\Bar{v}(t)}:=l_{v^i} (t,x,u), ~ \forall t \in [\tau_{i-1}, \tau_i), i=1,\cdots, N+1,
	\end{aligned}
  \label{Eq_lvi}
\end{equation}
and $x(\cdot)$ is a solution of (\ref{Eq_switching function}) under the control signal $u(\cdot)$ and the switching signal $\Bar{v}(\cdot)$, starting from $x(t_0)=x_0$.

The  SOCP is to find the optimal schedule  $\sigma^*=(q^*, \tau^*)$  that not only satisfies the boundary conditions but also minimizes the cost described in (\ref{Eq_Cost SW}).

\section{Embedded Optimal Control Formulation}
 
Most conventional numerical optimization algorithms based on Pontryagin's minimum principle are developed for problems with continuous decision variables. To take advantage of such algorithms, the SOCP is transformed into a continuous optimal control problem by embedding the switching variable into a continuous function $v:\mathbb{R}_{\geq 0}\rightarrow \mathcal{V}:=[0,1] $. The corresponding  embedded dynamics are expressed as
\begin{equation}
		\dot{x} =f(t,x, U ),
  \label{Eq_Em function}
\end{equation}
where $U:=[u_0^\top,u_1^\top, v]^\top$ denotes the augmented control, and $f(t,x,U) :=[1-v ]f_0(t,x,u_0) +v f_1(t,x,u_1) $.

The EOCP is then defined as (see \cite{SCC.Bengea.DeCarlo2005})

 \begin{equation*}
	\begin{aligned}
		\underset{U(\cdot)}{\min} ~~& J(t_0,x_0 ,t_f,U(\cdot))  \\
		\text{subject to}: & \\
		&  \text{(i) $(t_0,x_0,t_f,x_f) \in \mathcal{B}\subseteq \mathbb{R}^{2n+2} $}, \text{and}\\
		&  \text{(ii) $v(t)\in  \mathcal{V} , u_0(t), u_1(t) \in\mathcal{U} \in \mathbb{R}^m$,}\\
        &  ~ \forall t \in[t_0,t_f],
	\end{aligned}
\end{equation*} 
with the cost function
 \begin{equation}
	\begin{aligned}
 J(t_0,x_0 ,t_f,U(\cdot)):= &\int_{ t_0 }^{ t_f  } L(t,x(t),U(t))\mathrm{d}t\\
 &+ K(t_0,x_0,t_f,x_f),
	\end{aligned}
\end{equation}
where $L(t,x,U) =[1-v ]l_0(t,x,u_0) +v l_1(t,x,u_1)$, and $x(\cdot)$ is a solution of (\ref{Eq_Em function}) under the control signal $U(\cdot)$ starting from $x(t_0)=x_0$. 

Note that constraint (iii), corresponding to the dwell-time constraint in the SOCP, has been removed to render the EOCP solvable through Pontryagin’s minimum principle.

This paper focuses  on switching systems with two subsystems. However, the approach can be readily extended to switched systems with an arbitrary number of subsystems. An example of extending the embedded approach to a switching system with three subsystems can be found in \cite{SCC.Bengea.DeCarlo2005}.
\section{Modified Embedded Optimal Control Formulation}
While embedding enables the use of traditional solution techniques, optimal solutions of EOCPs can take any value in the interval of $[0,1]$, which cannot be directly interpreted as a solution for the SOCP. However, bang-bang solutions of the EOCP, where the embedded signal corresponds to one active mode at a time, are feasible for the SOCP.
In \cite{SCC.Abudia.Harlan.ea2020}, a MEOCP is developed where it is shown that adding a suitable auxiliary concave function to the cost function forces the EOCP to have a bang-bang optimal solution. In the  MEOCP formulation, the cost function is  modified as (see \cite{SCC.Abudia.Harlan.ea2020})
\begin{equation}
	\begin{aligned}
J(t_0  ,x_0 ,t_f,U(\cdot)):=& \int_{ t_0 }^{ t_f  }  \{L(t,x(t),U(t)) + L_{v} ( v ) \}\mathrm{d}t\\
& + K(t_0,x_0,t_f,x_f), 
	\end{aligned}
\end{equation}
where $  L_{v} ( v): \mathcal{V} \rightarrow \mathbb{R} $  is a concave function such that  $ L_{v } ( v )=0$ whenever $v \in \bar{\mathcal{V}} $. 

The Hamiltonian for the MEOCP is given by
\begin{equation}
	\begin{aligned}
H(t ,x ,p, U(\cdot))=& \langle~  p~,~ f(t,x,U)  ~\rangle + L(t,x,U)+ L_{v } ( v ),
	\end{aligned}
\end{equation}
where $p$ represents the costate. 

Based on Pontryagin's minimum principle \cite{SCC.Kopp1962}, the Hamiltonian is subject to the following necessary conditions along the optimal state trajectory $x^*(\cdot)$ and optimal costate trajectory $p^*(\cdot)$:
 \begin{subequations} 
    	\begin{align}
& H(x^*(t),p^*(t),U^*(t),t)\leq  H(x^*(t),p^*(t),U(t),t)\\
& \dot{p}^*(t)=-H_x(x^*(t),p(t),U^*(t),t)\\
& H(t_f)=0\\
& p(t_f)=0^n 
	\end{align}
 \end{subequations}

Since along the optimal trajectory, $U^*(t)=[{u_0^*}^{\top}(t),{u_1^*}^{\top}(t), {v ^*}(t)]^{\top}$  for all $t$ minimizes the Hamiltonian among all other controllers, we can conclude that the function $v \rightarrow H(x^*(t),p^*(t), [{u_0^*}^{\top}(t),{u_1^*}^{\top}(t), v]^{\top},t)$ is minimized by the optimal mode signal $ v ^*(t)$  for all $t$.

Define $U^*_{v}(t):=[{u_0^*}^\top(t),{u_1^*}^\top(t), v]^\top$. Since 
$ v \rightarrow H(t ,x ,p, U^*_{v},t)=\langle~  p~,~ f(t,x,U^*_{v},t)  ~\rangle + L(t,x,U^*_{v})+ L_{v } ( v )$ is a sum of an affine function $\langle~  p~,~ f(t,x,U^*_{v},t)  ~\rangle + L(t,x,U^*_{v})$  and concave function $L_{v } ( v )$, it is concave for all $t$.

Since a continuous concave function over a compact set achieves its minimum on  the boundary of the compact set (see  \cite{SCC.Zangwil1967}, Theorem 3), we can conclude that $v \rightarrow H(x^*(t),p^*(t), U^*_{v},t)$ achieves its minimum at the boundary of the range of the embedded mode signal. For switching systems with two sub-systems, the boundary of $[0,1]$ is $\{0,1\}=\bar{\mathcal{V}}$. As a result, the concave auxiliary function $L_v$ forces the MEOCP to have a bang-bang solution, which can be identified directly as a feasible solution of the SOCP. However, the solution does not guarantee satisfaction of the dwell time constraints (iii).
In  \cite{SCC.Abudia.Harlan.ea2020}  it has been demonstrated that dwell-time constraints could be heuristically met by adjusting the auxiliary cost function. However, this approach heavily depends on the numerical optimization algorithm and is therefore deemed unreliable. An alternative approach to meet the dwell time constraint has been adopted in the following section. The approach involves the introduction of an additional layer to post-process the solution of the MEOCP.

\section{Dwell-Time  Filtering}

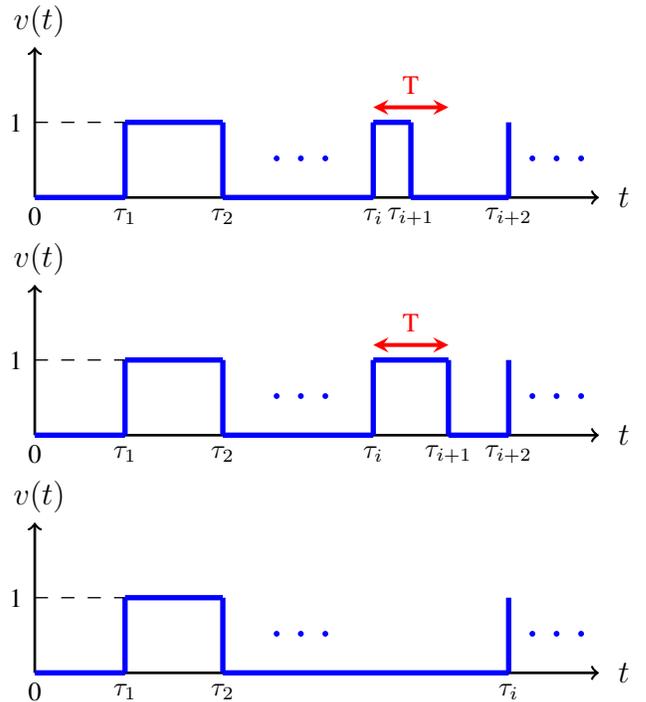
\begin{figure}[ht]
    \centering
\begin{tikzpicture}
    \draw[->, line width=1pt] (0,0) -- (7.5,0) node[right] {\fontsize{12}{14} $t$};
    \draw[->, line width=1pt] (0,0) -- (0,2) node[above] {\fontsize{12}{14} $v(t)$};
    
    \draw[line width=2pt,blue] (0,0) -- (1.2,0);
    \draw[line width=2pt,blue] (1.2,0) -- (1.2,1);
 
    \draw[line width=2pt,blue] (1.2,1) -- (2.5,1);
    \draw[line width=2pt,blue] (2.5,1) -- (2.5,0);
    
    \draw[line width=2pt,blue] (2.5,0) -- (4.5,0);
    \draw[line width=2pt,blue] (4.5,0) -- (4.5,1);

    \draw[line width=2pt,blue] (4.5,1) -- (5,1);
    \draw[line width=2pt,blue] (5,1) -- (5,0);

    \draw[line width=2pt,blue] (5,0) -- (6.3,0);
    \draw[line width=2pt,blue] (6.3,0) -- (6.3,1);

    \draw[>=stealth, <->, line width=1.5pt, red] (4.5,1.2) -- (5.5,1.2);
    \node[text=red, font=\fontsize{10}{14}\selectfont] at (5,1.5) {T};
    \draw[-, dash pattern=on 5pt off 5pt,line width=0.5pt ] (0,1) -- (1.2,1) ;

    \node[font=\fontsize{10}{14}\selectfont] at (0,-0.25) {0};
    \node[font=\fontsize{10}{14}\selectfont] at (-0.25, 1) {1};
    \node[font=\fontsize{10}{14}\selectfont] at (1.2,-0.25) {$\tau_1$};
    \node[font=\fontsize{10}{14}\selectfont] at (2.5,-0.25) {$\tau_2$};
    \node[font=\fontsize{10}{14}\selectfont] at (4.5,-0.25) {$\tau_i$};
    \node[font=\fontsize{10}{14}\selectfont] at (5,-0.25) {$\tau_{i+1}$};
    \node[font=\fontsize{10}{14}\selectfont] at (6.3,-0.25) {$\tau_{i+2}$};
    \node[font=\fontsize{20}{14}\selectfont, blue] at (3.6,0.5) {$\cdots$};
    \node[font=\fontsize{20}{14}\selectfont, blue] at (7,0.5) {$\cdots$};
\end{tikzpicture}
\begin{tikzpicture}
    \draw[->, line width=1pt] (0,0) -- (7.5,0) node[right] {\fontsize{12}{14} $t$};
    \draw[->, line width=1pt] (0,0) -- (0,2) node[above] {\fontsize{12}{14} $v(t)$};
    
    \draw[line width=2pt,blue] (0,0) -- (1.2,0);
    \draw[line width=2pt,blue] (1.2,0) -- (1.2,1);
    
    \draw[line width=2pt,blue] (1.2,1) -- (2.5,1);
    \draw[line width=2pt,blue] (2.5,1) -- (2.5,0);
    
    \draw[line width=2pt,blue] (2.5,0) -- (4.5,0);
    \draw[line width=2pt,blue] (4.5,0) -- (4.5,1);

    \draw[line width=2pt,blue] (4.5,1) -- (5.5,1);
    \draw[line width=2pt,blue] (5.5,1) -- (5.5,0);

    \draw[line width=2pt,blue] (5.5,0) -- (6.3,0);
    \draw[line width=2pt,blue] (6.3,0) -- (6.3,1);

    \draw[>=stealth, <->, line width=1.5pt, red] (4.5,1.2) -- (5.5,1.2);
    \node[text=red, font=\fontsize{10}{14}\selectfont] at (5,1.5) {T};
    \draw[-, dash pattern=on 5pt off 5pt,line width=0.5pt ] (0,1) -- (1.2,1) ;

    \node[font=\fontsize{10}{14}\selectfont] at (0,-0.25) {0};
    \node[font=\fontsize{10}{14}\selectfont] at (-0.25, 1) {1};
    \node[font=\fontsize{10}{14}\selectfont] at (1.2,-0.25) {$\tau_1$};
    \node[font=\fontsize{10}{14}\selectfont] at (2.5,-0.25) {$\tau_2$};
    \node[font=\fontsize{10}{14}\selectfont] at (4.5,-0.25) {$\tau_i$};
    \node[font=\fontsize{10}{14}\selectfont] at (5.5,-0.25) {$\tau_{i+1}$};
    \node[font=\fontsize{10}{14}\selectfont] at (6.3,-0.25) {$\tau_{i+2}$};
    \node[font=\fontsize{20}{14}\selectfont, blue] at (3.6,0.5) {$\cdots$};
    \node[font=\fontsize{20}{14}\selectfont, blue] at (7,0.5) {$\cdots$};   
\end{tikzpicture}
\begin{tikzpicture}
    \draw[->, line width=1pt] (0,0) -- (7.5,0) node[right] {\fontsize{12}{14} $t$};
    \draw[->, line width=1pt] (0,0) -- (0,2) node[above] {\fontsize{12}{14} $v(t)$};
    
    \draw[line width=2pt,blue] (0,0) -- (1.2,0);
    \draw[line width=2pt,blue] (1.2,0) -- (1.2,1);
    
    \draw[line width=2pt,blue] (1.2,1) -- (2.5,1);
    \draw[line width=2pt,blue] (2.5,1) -- (2.5,0);  

    \draw[line width=2pt,blue] (2.5,0) -- (6.3,0);
    \draw[line width=2pt,blue] (6.3,0) -- (6.3,1);

    \draw[-, dash pattern=on 5pt off 5pt,line width=0.5pt ] (0,1) -- (1.2,1) ;

    \node[font=\fontsize{10}{14}\selectfont] at (0,-0.25) {0};
    \node[font=\fontsize{10}{14}\selectfont] at (-0.25, 1) {1};
    \node[font=\fontsize{10}{14}\selectfont] at (1.2,-0.25) {$\tau_1$};
    \node[font=\fontsize{10}{14}\selectfont] at (2.5,-0.25) {$\tau_2$};
    \node[font=\fontsize{10}{14}\selectfont] at (6.3,-0.25) {$\tau_i$};
    \node[font=\fontsize{20}{14}\selectfont, blue] at (3.6,0.5) {$\cdots$};
    \node[font=\fontsize{20}{14}\selectfont, blue] at (7,0.5) {$\cdots$};  
\end{tikzpicture}

\caption{The typical switching mode signal in order from top to bottom: original switching mode
signal, replace mode $1$ during interval $[\tau_i,\tau_i+T) $, replace mode $0$ during interval $[\tau_i,\tau_i+T) $}
\label{Fig. Mode removal}
\end{figure}

To impose the dwell-time constraint, a filtering layer is added to the solution of the  MEOCP. This layer ensures compliance with the dwell-time constraint by eliminating any switching times that violate it. An illustration of this process is provided in Fig.~\ref{Fig. Mode removal}. It is evident from the figure that the switching mode signal $v(\cdot)$ fails to satisfy the dwell-time constraint at time $\tau_{i+1}$,
as indicated by the fact that  $\tau_{i+1}-\tau_i < T$.
Consequently, the switch at $\tau_{i+1}$ needs to be removed in such a way that the resulting increase in the total cost is as small as possible.
Therefore, the primary inquiry revolves around determining which mode should be selected during the interval $[\tau_i, \tau_i+T)$.

If there were a single dwell-time violation, the question above could be answered by simply computing the total cost for each mode. However, if there are multiple violations, then the number of mode combinations that need to be evaluated grow geometrically with respect to the number of violations. As a result, an approach based on computation of total cost is not numerically feasible in general. In this paper, inspired by the mode insertion method developed in \cite{SCC.Axelsson.Wardi.ea2007}  and \cite{SCC.Wardi.Egerstedt.ea2014}, the mode selection question is answered using the gradient of cost with respect to the insertion time.

To illustrate the mode insertion method consider $\sigma=( q, \tau)$ in Definition 2 as a feasible schedule for the system  in (\ref{Eq_switching function})  and express it as
\begin{equation}
	\begin{aligned}
\dot{x}=f_{v^i} (t,x,u), ~~ \forall t \in [\tau_{i-1}, \tau_i), i=1,\cdots, N+1
	\end{aligned}
 \label{Eq_fvi}
\end{equation}

By defining $F_R$ and $L_R$ as right-hand sides (RHS) of  equations (\ref{Eq_fvi}) and (\ref{Eq_lvi}) respectively, the costate $p(t)$ can be shown to be a solution of the differential equation
\begin{equation}
	\begin{aligned}
\dot{p}=-\left( \pdv{F_R}{x} (x,u,t)\right)^\top p-\left(\pdv{L_R}{x} (x,u,t)\right)^{\top},
	\end{aligned}
  \label{Eq_costate}
\end{equation}
with the boundary condition $p(t_f)=0$~\cite{SCC.Wardi.Egerstedt.ea2014}.
Suppose that we have a feasible schedule $\sigma=( q, \tau)$  and we want to insert the new mode $ \alpha \in \bar{\mathcal{V}}$  during the time interval $[ s,s+\lambda)$ for selected  $s \in [t_0,t_f)$ and $\lambda >0 $. By fixing insertion time $s$ and insertion mode $\alpha$, the cost in (\ref{Eq_Cost SW}) can be seen as a function of $\lambda$,  denoted by $ J_{\sigma,s,\alpha}(\lambda) $.  
Denote the right-hand side derivative of the $ J_{\sigma,s,\alpha}(\lambda) $ at $\lambda=0$ by $D_{\sigma,s,\alpha}$, which can be calculated as (see \cite{SCC.Axelsson.Wardi.ea2005})
\begin{equation}
\begin{aligned}
D_{\sigma,s,\alpha}:=&\pdv{J_{\sigma,s,\alpha}}{\lambda^+} (0) \\
=&p(s)^\top \left \{ f_\alpha (s, x(s),u(s))- f_{v(s)}(s, x(s),u(s)) \right \}. 
\end{aligned}
\end{equation}

The derivative $D_{\sigma,s,\alpha}$ is called the insertion gradient.

To answer the primary question of which mode needs to be selected during the interval $[\tau_i, \tau_i+T)$, as illustrated in Fig.~\ref{Fig. Mode removal}, we resolve the optimization problem

          \begin{equation}
	      \begin{aligned}
              \alpha:= \underset{v \in \bar{\mathcal{V}}}{\min} ~~& \int_{\tau_i}^{\tau_i+T}D_{\sigma,t,v}~ \mathrm{d}t,
	      \end{aligned}
       \label{Eq_best mode}
         \end{equation}
where $\alpha$ denotes the switching mode than needs to be selected in the interval $[\tau_i, \tau_i+T)$.

The mode that minimizes the integral in equation~(\ref{Eq_best mode}) is used to replace the  modes in the interval $[\tau_i, \tau_i+T)$  based on the intuition that choosing the mode with the smallest average derivative over the interval $[\tau_i, \tau_i+T)$, will minimize the impact of the change on the optimal cost.

To calculate $\alpha$ numerically, one can choose an appropriate number of sub-intervals within the interval $[\tau_i, \tau_i+T)$ to determine the step size $\mathrm{d}t$. Then, for each mode, calculate the insertion gradient in each sub-interval and apply the trapezoidal rule to compute the integral in (\ref{Eq_best mode}). Finally, select the mode that minimizes the integral. 

Consider $\sigma^*=( q^*, \tau ^*)$ as a solution of the MEOCP. The filter designed to ensure compliance with the dwell-time constraint can be implemented through the algorithm outlined in Algorithm \ref{alg_1} to achieve filtered schedule $\sigma_f$.
Initially, the embedded schedule undergoes a check for any violations of dwell-time constraints. If no violations exist, the filtered schedule $\sigma_f$ remains identical to the embedded schedule $\sigma^*$. However, if violations are detected, the first instance of a dwell-time constraint being breached is pinpointed and denoted as $\tau_{i+1}$, identified by the condition $\tau_{i+1}-\tau_i < T$. Subsequently, the switch at time $\tau_{i+1}$ is eliminated, and mode $\alpha$ is determined from (\ref{Eq_best mode}) for the interval $[\tau_i,\tau_i+T)$ (see Fig.~\ref{Fig. Mode removal}). To ensure consecutive sequencing, the switching times are renumbered as $\tau_1, \tau_2, \cdots$. This process is iterated for the subsequent steps until no dwell-time constraint violations are detected.
\begin{algorithm}
\caption{ Filtering the solution of MEOCP }
\label{alg_1}
\begin{algorithmic}[1] 
    \STATE Set $k=1$
    \STATE Set $\sigma=\sigma^*$
        \IF{$\sigma$ does not satisfy the dwell-time constraint}
            \STATE Compute $j:=\min \{ i\geq k \mid  \tau_{i+1}-\tau_i < T \}$
            \STATE Removes all the switching times of $\sigma$ that is occurring between $\tau_j$ and $\tau_j+T$ 
            \STATE Compute $\alpha$:
          \begin{equation}
	      \begin{aligned}
              \alpha:= \underset{v \in \bar{\mathcal{V}}}{\min} ~~& \int _{\tau_j}^{\tau_j+T} D_{\sigma,t,v} \mathrm{d}t
	      \end{aligned}
         \end{equation}
         \STATE Set $\alpha$ as a mode of $\sigma$ during the interval $[t_j,t_j+T)$ 
         \STATE Renumbering the switching times  as $\tau_1, \tau_2,\cdots$  to make them consecutive.
         \STATE Set $ k=j+1$ and go to step 3
        \ENDIF
        \STATE $\sigma_{f}=\sigma$
\end{algorithmic}
\end{algorithm}
In contrast to the approach outlined in \cite{SCC.Abudia.Harlan.ea2020}, which attempts to meet the dwell-time constraint by adjusting constants within the auxiliary cost function heuristically, the filter utilized in this study offers a straightforward approach that delivers a schedule that, by construction, satisfies the dwell-time constraint.
 It should be noted that in spite of the filter,  mode deletion and replacement distorts the trajectory of the system. As a consequence, the post-replacement segment of the trajectory is no longer optimal.
To mitigate such loss of optimality, if computational resources allow, after completing step 7 in Algorithm \ref{alg_1}, one can re-solve the MEOCP over the interval $[t_j+T, t_f]$, with the initial condition set as $x_0 = x(t_j+T)$, renumber the switching times, and repeat the process after every mode replacement.

\section{Simulation results}
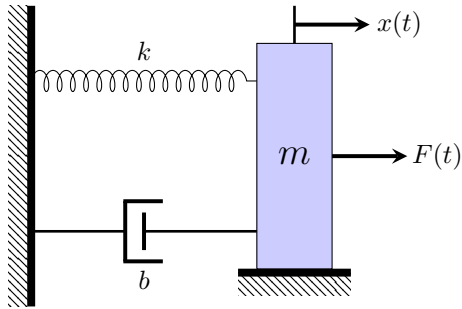
\begin{figure}[ht]
\centering
\begin{tikzpicture}

    \fill[pattern=north west lines] (-0.3,-2) rectangle (0,2); 
    \draw[line width=3pt] (0,-2) -- (0,2);   
    \fill[pattern=north west lines] (2.75,-1.85) rectangle (4.25,-1.55); 
    \draw[line width=3pt] (2.75,-1.55) -- (4.25,-1.55);   

    \filldraw[fill=blue!20] (3,-1.5) rectangle (4,1.5) node[ text=black, font=\fontsize{15}{14}\selectfont ] at (3.5,0) {$m$} ;
    
    \draw[decorate,decoration={coil,segment length=5pt,amplitude=4pt}] (0,1) -- (3,1) node[text=black, font=\fontsize{10}{14}\selectfont]  at ( 1.5, 1.4 ) {$k$};

\draw[line width=1pt] (0,-1) -- (1.25,-1);  
\draw[line width=1.5pt] (1.5,-1.25) -- (1.5,-0.75);   
\draw[line width=1.5pt] (1.25,-1.4) -- (1.25,-0.6);   
\draw[line width=1.5pt] (1.225,-0.6) -- (1.75,-0.6);   
\draw[line width=1.5pt] (1.225,-1.40) -- (1.75,-1.40); 
\draw[line width=1pt] (1.5,-1) -- (3,-1) node[text=black, font=\fontsize{10}{14}\selectfont]  at ( 1.5, -1.65 ) {$b$}; 

    \draw[>=stealth, ->,black, line width=1.5pt] (4,0) -- (5,0) node[text=black, font=\fontsize{10}{14}\selectfont]  at ( 5.4, 0 ) {$F(t)$};

    \draw[line width=1pt] (3.5,1.5) -- (3.5,2);
    \draw[>=stealth, ->,black, line width=1.5pt] (3.5,1.75) -- (4.5,1.75)  node[text=black, font=\fontsize{10}{14}\selectfont]  at ( 4.9, 1.75 ) {$x(t)$};
    
\end{tikzpicture}
	\caption{Mass-Spring-Damper System}
	\label{Fig. MDS}
\end{figure}

The mass-spring-damper system depicted in Fig.~\ref{Fig. MDS} has been employed for simulation purposes.
The parameters of the system are $m=1 kg$, $k=0.1 N/m$, and $b=0.1 kg/s$. The objective is to design a controller that achieves a mass position of $1m$ with zero velocity under the constraint that the force can only take values of $ -0.2 N$  or $ 0.2 N$.
The equation of motion of the mass-spring-damper system is 
 \begin{equation}
	\begin{aligned}
m \ddot{x}+b \dot{x}+kx=F,
	\end{aligned}
\end{equation} 
written as a  linear time-invariant (LTI) system
\begin{equation}
\dot {x}(t)
=
  \begin{pmatrix}
0 &  1 &  \\
-k/m & -b/m  \\
\end{pmatrix}
x(t)+
  \begin{pmatrix}
0 \\
1/m  \\
\end{pmatrix}
u(t),
\end{equation}
The force constraint makes the system a switching control system with two modes corresponding to $u(t)=F \in \{-0.2,0.2\}$. The initial condition is selected as $x=[0,0]^{\top}$. The cost and auxiliary cost functions are selected as $L=a\{(x_1-1)^2-x_2^2\}$ and $L_v=b(v-v^2)$ with $a=4$ and $b=1$.

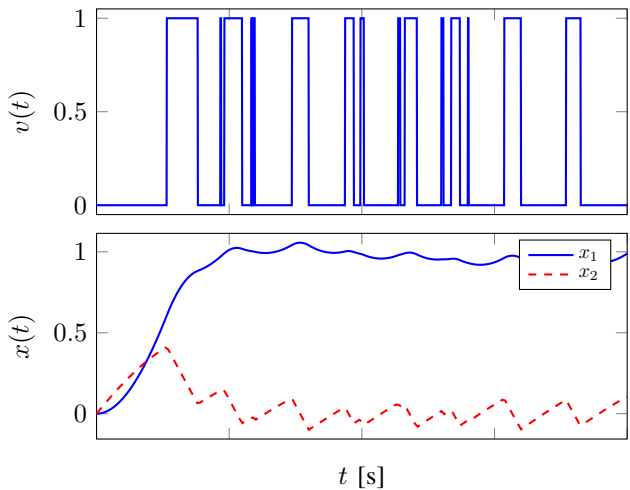
\begin{figure}[ht]
\hspace*{3em}\begin{tikzpicture}[trim axis left]
    \begin{axis}[
        xticklabels=\empty,
        ylabel={$v(t)$},
        legend pos = north east,
        legend style={nodes={scale=0.75, transform shape}},
        enlarge y limits=0.05,
        enlarge x limits=0,
        width=\linewidth,
        height=0.5\linewidth,
        ]
        \pgfplotsinvokeforeach {1} {\addplot+ [thick, mark=none] table [x index=0, y index=#1] {Data/Figs/V_emb.dat};}
        \legend{}
            \end{axis}
\end{tikzpicture}
\hspace*{3em}\begin{tikzpicture}[trim axis left]
    \begin{axis}[
        xlabel={$t$ [s]},
         xticklabels=\empty,
        ylabel={$x(t)$},
        legend pos = north east,
        legend style={nodes={scale=0.75, transform shape}},
        enlarge y limits=0.05,
        enlarge x limits=0,
        width=\linewidth,
        height=0.5\linewidth,
        ]
        \addplot+ [thick, mark=none, blue] table [x index=0, y index=1] {Data/Figs/X_emb.dat};
        \addplot+ [thick, mark=none, red, dashed] table [x index=0, y index=2] {Data/Figs/X_emb.dat};
        \legend{$x_1$, $x_2$}
    \end{axis}
\end{tikzpicture}
	\caption{  The solution of MEOCP (without a dwell-time constraint). The plot on top shows the switching mode signal as a function of time and the  plot at the bottom shows the states of the systems as a function of time.}
	\label{Fig.MEOCP_solution}
\end{figure}

\begin{figure}[ht]
\hspace*{3em}\begin{tikzpicture}[trim axis left]
    \begin{axis}[
         xticklabels=\empty,
        ylabel={$L(t)$},
        legend pos = north east,
        legend style={nodes={scale=0.75, transform shape}},
        enlarge y limits=0.05,
        enlarge x limits=0,
        width=\linewidth,
        height=0.5\linewidth,
        ]
        \addplot+ [thick, mark=none, blue] table [x index=0, y index=1] {Data/Figs/L_cost.dat};
        \addplot+ [thick, mark=none, red, dashed] table [x index=0, y index=2] {Data/Figs/L_cost.dat};
        \addplot+ [thick, mark=none, black, dash pattern=on 3pt off 2pt on 1pt off 2pt] table [x index=0, y index=3] {Data/Figs/L_cost.dat};
        \legend{$T=0s$, $T=0.1s$, $T=0.2s$}
    \end{axis}
\end{tikzpicture}
\hspace*{3em}\begin{tikzpicture}[trim axis left]
    \begin{axis}[
        xlabel={$t$ [s]},
        ylabel={$J(t)$},
        legend pos = south east,
        legend style={nodes={scale=0.75, transform shape}},
        enlarge y limits=0.05,
        enlarge x limits=0,
        width=\linewidth,
        height=0.5\linewidth,
        ]
        \addplot+ [thick, mark=none, blue] table [x index=0, y index=1] {Data/Figs/Cost.dat};
        \addplot+ [thick, mark=none, red, dashed] table [x index=0, y index=2] {Data/Figs/Cost.dat};
        \addplot+ [thick, mark=none, black, dash pattern=on 3pt off 2pt on 1pt off 2pt] table [x index=0, y index=3] {Data/Figs/Cost.dat};
        \legend{$T=0s$, $T=0.1s$, $T=0.2s$}
    \end{axis}
\end{tikzpicture}
	\caption{Cost for the solution of the MEOCP without a dwell-time filter, i.e., $T=0 s$,  and with a dwell-time filter with a $T=0.1 s$ and $T=0.2 s$. The plot on top shows the Lagrangian (running cost) as a function of time and the  plot at the bottom shows the cost as a function of time.}
	\label{Fig.Cost}
\end{figure}
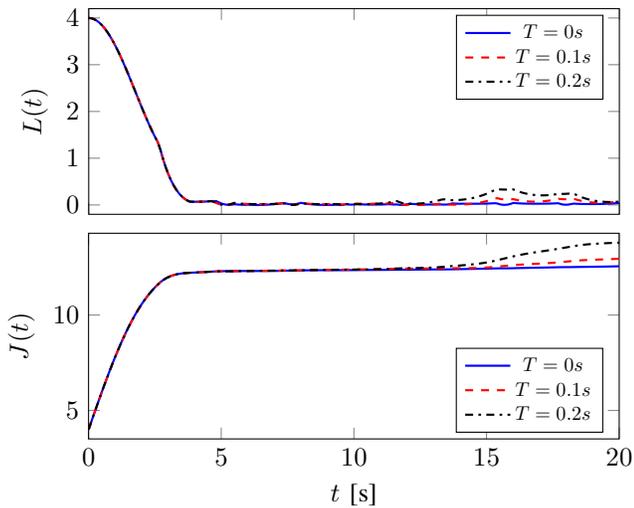

For solving the MEOCP, the GPOPS-II software is employed, which utilizes the 
adaptive Gaussian quadrature collocation method for solving multiple-phase optimal control problems \cite{SCC.Patterson.Rao2014}.
 The solution of the MEOCP is shown in  Fig.~\ref{Fig.MEOCP_solution}.  Note that the solution includes a high-frequency switching signal that may not be achievable in many practical applications due to a dwell time constraint. The effectiveness of the adopted filter is tested for a dwell time constraint with $T=0.1 s$ and $T=0.2 s$. 
The resulting switching signals and the corresponding state trajectories are shown in Fig. \ref{Fig.MEOCP+f0.1} and Fig.~\ref{Fig.MEOCP+f0.2}, respectively.
Fig.~\ref{Fig.Cost} compares the total cost for different dwell-time constraints.
 %

\begin{figure}
\hspace*{3em}\begin{tikzpicture}[trim axis left]
    \begin{axis}[
        xticklabels=\empty,
        ylabel={$v(t)$},
        legend pos = north east,
        legend style={nodes={scale=0.75, transform shape}},
        enlarge y limits=0.05,
        enlarge x limits=0,
        width=\linewidth,
        height=0.5\linewidth,
        ]
        \pgfplotsinvokeforeach {1} {\addplot+ [thick, mark=none] table [x index=0, y index=#1] {Data/Figs/V_filter01.dat};}
        \legend{}
    \end{axis}
\end{tikzpicture}
\hspace*{3em}\begin{tikzpicture}[trim axis left]
    \begin{axis}[
        xlabel={$t$ [s]},
         xticklabels=\empty,
        ylabel={$x(t)$},
        legend pos = north east,
        legend style={nodes={scale=0.75, transform shape}},
        enlarge y limits=0.05,
        enlarge x limits=0,
        width=\linewidth,
        height=0.5\linewidth,
        ]
       \addplot+ [thick, mark=none, blue] table [x index=0, y index=1] {Data/Figs/X_filter01.dat};
        \addplot+ [thick, mark=none, red, dashed] table [x index=0, y index=2] {Data/Figs/X_filter01.dat};
        \legend{$x_1$, $x_2$}
    \end{axis}
\end{tikzpicture}
	\caption{The results after adding the filter (for a dwell-time constraint $T=0.1 s$) to the solution of  MEOCP. The plot on top shows the switching mode signal as a function of time and the  plot at the bottom shows the states of the systems as a function of time.}
	\label{Fig.MEOCP+f0.1}
\end{figure}
\begin{figure}
\hspace*{3em}\begin{tikzpicture}[trim axis left]
    \begin{axis}[
        xticklabels=\empty,
        ylabel={$v(t)$},
        legend pos = north east,
        legend style={nodes={scale=0.75, transform shape}},
        enlarge y limits=0.05,
        enlarge x limits=0,
        width=\linewidth,
        height=0.5\linewidth,
        ]
        \pgfplotsinvokeforeach {1} {\addplot+ [thick, mark=none] table [x index=0, y index=#1] {Data/Figs/V_filter02.dat};}
        \legend{}
            \end{axis}
\end{tikzpicture}
\hspace*{3em}\begin{tikzpicture}[trim axis left]
    \begin{axis}[
        xlabel={$t$ [s]},
         xticklabels=\empty,
        ylabel={$x(t)$},
        legend pos = north east,
        legend style={nodes={scale=0.75, transform shape}},
        enlarge y limits=0.05,
        enlarge x limits=0,
        width=\linewidth,
        height=0.5\linewidth,
        ]
        \addplot+ [thick, mark=none, blue] table [x index=0, y index=1] {Data/Figs/X_filter02.dat};
        \addplot+ [thick, mark=none, red, dashed] table [x index=0, y index=2] {Data/Figs/X_filter02.dat};
        \legend{$x_1$, $x_2$}
    \end{axis}
\end{tikzpicture}
	\caption{The results after adding the filter (for a dwell-time constraint $T=0.2 s$) to the solution of  MEOCP. The plot on top shows the switching mode signal as a function of time and the  plot at the bottom shows the states of the systems as a function of time. }
	\label{Fig.MEOCP+f0.2}
\end{figure}
\section{Discussion}
Imposing a stricter dwell-time constraint naturally leads to fewer switching occurrences, as shown in Fig.\ref{Fig.MEOCP+f0.1} and Fig.\ref{Fig.MEOCP+f0.2}, along with diminished performance, as reflected by the increased costs in Fig.~\ref{Fig.Cost}. This connection between the dwell-time constraint and the resulting cost emphasizes the trade-off between utilizing high-frequency actuators and ensuring the performance of the control process. 

It is important to note that the auxiliary cost does not affect the total cost in a true bang-bang solution because $L_v(v) = 0$ whenever $v \in \bar{\mathcal{V}}$. However, numerical approximation of a bang-bang solution can result in a solution that contains intermediate values of the switching signals. The presence of these intermediate values leads to a nonzero contribution to the auxiliary cost function, which affects the total cost. 
To apply the dwell-time filter, the intermediate values in the switching signals of MEOCP are first rounded to the nearest mode numbers which in this paper is $0$ or $1$.

The method developed in~\cite{SCC.Abudia.Harlan.ea2020} aims to satisfy the dwell-time constraint by adjusting the coefficients of an auxiliary cost function. Such an approach can increase the average dwell time but may still result in isolated dwell-time violations. Furthermore, this approach relies solely on the intermediate values artificially introduced by the solver to adjust the dwell time, and as such, is not reliable. In contrast, the post-filter method adopted in this paper avoids such violations.
\section{Conclusion}
This paper presents a post-filtered embedding approach to solving SOCPs with dwell time constraints. First, the switching mode signal is incorporated into a broader domain by substituting it with a suitable continuous function.
Then, by adding the auxiliary cost function, the EOCP is forced to have a bang-bang solution. Although the solution of MEOCP can be identified directly as a solution of SOCP, it will not guarantee the dwell time constraint, which is crucial in practical applications. Unlike previous approaches that adjust the coefficient of the auxiliary function heuristically, which is unreliable and highly dependent on the optimization algorithm, a mode-insertion-inspired filter layer is added to the solution of MEOCP, in this method resulting in a switching mode signal that satisfies the dwell-time constraint. 

The simulation results show that the approach is effective, but  only for small dwell time constraints; otherwise, the solution ends up with very few switches, leading to poor performance. Another drawback of the method is that after removing the switching instances and replacing them with a proper switching mode, the resulting solution is no longer optimal, even if the MEOCP is resolved for the remaining trajectory.

The  technique presented in this paper focuses on switching systems with only two subsystems. Extending this approach to systems with multiple subsystems, by an appropriate  defining an auxiliary concave function in higher dimensions, remains an area for future research. Another promising direction for future work is developing a post-filter that ensures compliance with an average dwell-time constraint. Additionally, exploring hybrid optimization algorithms that handle both continuous and discrete variables while incorporating dwell-time constraints offers an exciting avenue for further investigation.



\small
\bibliographystyle{IEEETrans.bst}
\bibliography{Ref/MyRef.bib}

\end{document}